\documentclass[11pt]{article}
\usepackage{amssymb,epsfig}
\newcommand{\2}{\vspace{0.5 cm}}

\newcommand{\qed}{\hfill$\Box$}
\newcommand{\pf}{{\bf Proof: }}
\newtheorem{theorem}{Theorem}[section]

\newtheorem{proposition}[theorem]{Proposition}
\newtheorem{problem}[theorem]{Problem}

\newtheorem{corollary}[theorem]{Corollary}
\newtheorem{conjecture}[theorem]{Conjecture}

\newcommand{\beq}{\begin{equation}}
\newcommand{\eeq}{\end{equation}}

\newcommand{\al}{\alpha}
\newcommand{\be}{\beta}

\begin{document}

\title{Mediated Digraphs and Quantum Nonlocality\footnote{Research of Gutin, Rafiey and Yeo was
supported in part by the Leverhulme Trust. Research of Gutin and
Rafiey was supported in part by the IST Programme of the European
Community, under the PASCAL Network of Excellence,
IST-2002-506778.}}

\date{}

\author{
G. Gutin\thanks{Department of Computer Science, Royal Holloway,
University of London, Egham, Surrey, TW20 0EX, UK,
gutin@cs.rhul.ac.uk}\and N. Jones \thanks{ Department of
Mathematics, Bristol University, University Walk, Bristol, BS8
1TW, UK, n.s.jones@bristol.ac.uk} \and A. Rafiey\thanks{Same as
Gutin apart from arash@cs.rhul.ac.uk} \and S.
Severini\thanks{Department of Mathematics and Department of
Computer Science, University of York, YO10 5DD, UK,
ss54@york.ac.uk} \and A. Yeo\thanks{Same as Gutin apart from
anders@cs.rhul.ac.uk}}

\maketitle

\begin{abstract}
A digraph $D=(V,A)$ is mediated if for each pair $x,y$ of distinct
vertices of $D$, either $xy\in A$ or $yx\in A$ or there is a vertex
$z$ such that both $xz,yz\in A.$ For a digraph $D$, $\Delta^-(D)$ is
the maximum in-degree of a vertex in $D$. The $n$th mediation number
$\mu (n)$ is the minimum of $\Delta^-(D)$ over all mediated digraphs
on $n$ vertices. Mediated digraphs and $\mu(n)$ are of interest in
the study of quantum nonlocality.

We obtain a lower bound $f(n)$ for $\mu(n)$ and determine infinite
sequences of values of $n$ for which $\mu(n)=f(n)$ and
$\mu(n)>f(n)$, respectively. We derive upper bounds for $\mu(n)$ and
prove that $\mu(n)=f(n)(1+o(1))$. We conjecture that there is a
constant $c$ such that $\mu(n)\le f(n)+c.$ Methods and results of
design theory and number theory are used.

\end{abstract}

\section{Introduction}
The class of mediated digraphs defined later in this section was
introduced in \cite{Jo} as a model in quantum mechanics. We define
and study an extremal parameter of digraphs in this class, the
$n$th mediation number. The parameter is of interest in the study
of quantum nonlocality.

The vertex (arc) set of a digraph $D$ will be denoted by $V(D)$
($A(D)$). For a digraph $D$ and $x\neq y\in V(D)$, we say that $x$
{\em dominates} $y$ if $xy\in A(D).$ All vertices that dominate
$x$ are {\em in-neighbors} of $x$; the set of in-neighbors is
denoted by $N^-(x).$ The number of in-neighbors of $x$ is the {\em
in-degree} of $x$. The {\em closed in-neighborhood} $N^-[x]$ is
defined as follows: $N^-[x]=\{x\}\cup N^-(x).$ We denote the
maximum in-degree of a vertex of a digraph $D$ by $\Delta^-(D).$
For standard terminology and notation on digraphs, see, e.g.,
\cite{bang2000}.

\begin{figure}[h]
\label{f}
\begin{center}
\epsfig{file=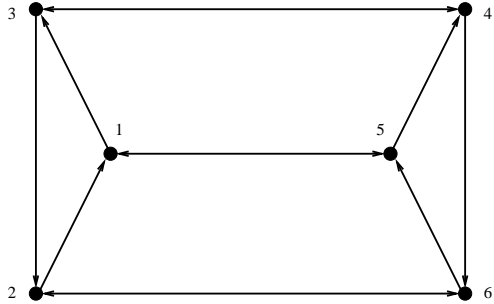, height=4cm}
\end{center}
\caption{A mediated digraph $H$ of order 6}
\end{figure}

A digraph $D$ is {\em mediated} if for every pair $x,y$ of vertices
there is a vertex $z$ such that both $x,y\in N^-[z]$ (possibly $z=x$
or $y$). Tournaments, doubly regular digraphs \cite{ioninJCT101} and
symmetric digraphs of diameter 2 are special families of mediated
digraphs. Fig. \ref{f} is an example of a mediated digraph.

The $n$th {\em mediation number} $\mu (n)$ is the minimum of
$\Delta^-(D)$ over all mediated digraphs on $n$ vertices. This
parameter is of interest in quantum mechanics as explained in the
next section.

The rest of the paper is organized as follows. Section \ref{qmsec}
provides a motivation for the study of mediated digraphs and the
$n$th mediation number. (One is not required to read the section
in order to understand the rest of the paper.) In Section
\ref{lower}, we obtain a lower bound $f(n)$ for $\mu(n)$, which is
proved to be sharp in the next two sections. Section
\ref{intersecting} is devoted to a characterization of $\mu(n)$ as
an extremal parameter of special families of sets. This allows us
to use some results from design theory. Section \ref{upper}
provides upper bounds for $\mu(n)$. We prove that
$\mu(n)=f(n)(1+o(1)),$ which is the central result of the paper
and is of importance for quantum nonlocality (see Section
\ref{qmsec}). In Section \ref{>sec} we show that $\mu(n)>f(n)$ for
an infinite number of values of $n.$   We conjecture that, in
fact, there is a constant $c$ such that $\mu(n)\le f(n)+c$ for
each $n\ge 1$ and pose the problem of checking whether $\mu(n)$ is
a monotonically increasing function.

\section{Mediated Digraphs in Quantum Mechanics}\label{qmsec}

Nonlocality is a fundamental, and curious, feature of quantum theory
which confused Einstein and continues to yield exciting results in
physics (there are numerous popular explanations of nonlocality, a
more technical review is \cite{Werner}  and in \cite{Bell} one can
find some of the crucial early papers in the field). The study of
nonlocality is sometimes helped by considering classical analogies:
it was in this endeavor that mediated digraphs were discovered (see
\cite{Jo} - in that paper mediated digraphs are called Totally
Paired Graphs). Consider two objects which are connected and then
suddenly sent to such widely separated locations that they can no
longer influence each other on relevant time scales. The results of
local measurements on each member of a pair of classical objects,
which have been connected and separated in this fashion, can be
correlated, depending on their relationship when they were together.
Perhaps, when they were together, the objects exchanged some
information, like a string of bits. The rough edges of a sheet of
paper torn in two remain correlated when the pieces are sent to
distantly separated locations: local measurements that are made on
them are connected. The correlations between the results of certain
sets of local measurements on some pairs of quantum objects, which
have been connected and separated, cannot be explained by allowing
only an exchange of a bit string when they started together. If one
studies the probability distributions of the different possible
outcomes of local measurements on sets of quantum objects, for
different local measurement settings, one cannot explain them by
common strings of information shared between the objects. This is an
aspect of nonlocality. Since it is the case that measurements on
quantum objects can show classical correlations but the reverse is
not true, there is a sense in which quantum objects have
correlations beyond allowed classical ones.

Nonlocality has been well studied for pairs of quantum objects but
less work has been undertaken for more than two \cite{Col,Seev}.
Let $i$ be some object which can be measured in one of two ways
but not both at once. For example, suppose we have an apparatus
which measures either the height or the width of $i$, but not the
two together. Let $x_i\in\{0,1\}$ be the measurements of $i$  and
let $a_i\in\{0,1\}$ be possible outcomes of this measurement. It
is natural to hold that the way that an object is measured affects
the result of a measurement: $a_i=a_i(x_i)$ (an object with height
`0' and width `1' will yield a result `1' when its width is
measured and `0' otherwise). If the measurement events occur at
space-like separated locations then the way that another object,
$j$, is measured elsewhere, $x_j$, cannot affect $a_i$ unless the
objects are exchanging information faster than light:
$a_i=a_i(x_i,x_j)=a_i(x_i,(x_j+1) \mbox{ mod } 2)$ ($a_i$ is
unaltered for any value of $x_j$).

A standard classical analogy for quantum nonlocality is as follows
(see \cite{Toner} and references therein). Classical separated
objects are allowed to cheat and exchange information, faster than
light, about the way they are to be measured. In this case the
outcomes of a measurement on the object $i$ could indeed depend on
the way the object $j$ is measured: $a_i=a_i(x_i,x_j)$. The
correlations present in sets of quantum objects can now be
classically approximated. One can ask how many bits of information
have to be exchanged between classical objects in order to fool an
experimentalist into thinking he/she is measuring a quantum state.
The classical objects are given an extra property; their
characteristics can depend on the way other, distant, objects are
measured. How much of this freedom need one allow in order to
produce scenarios which can have the same measurement results as
measurements of quantum objects?

In the scenario considered in \cite{Jo} (motivated by the
structure of probability spaces) each object (a vertex) knows how
it is to be measured and can send this information to other
vertices (an arc from source vertex to target vertex). This
information stays put on receipt and does not propagate around the
graph:  a vertex can only know the way another vertex is to be
measured by receiving an arc directly from that vertex (not via a
third party). The measurement results of vertices, given that
their properties might now depend on the way their neighbors are
to be measured, are more general and now have a chance to
reproduce quantum correlations. It was shown that it is necessary
that the vertices be connected as a mediated digraph if they are
to fool an experimentalist into thinking he/she is measuring a
quantum state. Within this model, a certain topology of
communication is a necessary classical property  in order to
simulate sets of quantum objects classically. Note that
sufficiency was not shown and this is now being studied (these
digraphs are patterns of faster than light communication - however
this violation of causality is not necessary and can be removed by
a randomization procedure \cite{Jo} Theorem 3).

Given that $n$ classical objects connected as any mediated digraph
can sometimes be at least as nonlocal as $n$ quantum mechanical
objects, it is interesting to find out how 'connected' these
digraphs are. If the digraphs are good analogues of quantum
nonlocality, then their structure should inform us about quantum
correlations. One would like to consider the least connected members
of the set of mediated digraphs - the least connected digraphs that
can still be at least as nonlocal as quantum states. In order to
achieve this one must have a good measure of connectivity: we
consider $\Delta^-(D).$ If an $n$ vertex digraph contains a vertex
which depends on the settings of lots of other vertices,
$\Delta^-(D)$ will be large: this defines a highly nonlocal pattern
- one vertex is highly correlated with many others. If all vertices
in a digraph are only connected to a few others, $\Delta^-(D)$ will
be small: such digraphs seem to have a form of short range
nonlocality. Proving that, for any $n$, there are mediated digraphs
which have $\Delta^-(D)$ scaling with $\sqrt{n}$ (Theorem
\ref{o(1)th} of this paper), shows that each object need only be
connected to a fraction of the set of objects which {\em diminishes}
as $n$ increases (as $1/\sqrt{n}$). As $n$ increases there exist
mediated digraphs in which each vertex becomes increasingly
localized with respect to the whole - this must be telling us
something about quantum nonlocality.

\section{Lower Bound for $\mu(n)$}\label{lower}

For a real $x$, let $\lceil x \rceil$ denote the least integer not
smaller than $x.$ Let $f(n)=\lceil
\frac{1}{2}(\sqrt{4n-3}-1)\rceil.$ The following proposition gives a
lower bound for $\mu(n)$, which is the exact value of $\mu(n)$ for
infinitely many values of $n$ (see Corollaries \ref{muproj} and
\ref{projth}).

\begin{proposition}\label{lbprop}
For each $n\ge 1$, we have $\mu(n)\ge f(n)$.
\end{proposition}
\pf Let $D$ be a mediated digraph and let $d=\Delta^-(D).$ If $D$
has just one vertex, the bound holds, so we may assume that $n\ge
2.$ By the definition of a mediated digraph, each pair $x,y$ of
vertices of $D$ belongs to the closed in-neighborhood of some
vertex. Let $d^-_1, d^-_2, \ldots , d^-_n$ be the in-degrees of
vertices $v_1,v_2,\ldots,v_n$ of $D$. Since a vertex $v_i$ has
${d^-_i \choose 2}+ d^-_i$ pairs of vertices in its closed
in-neighborhood and since $D$ has overall $n \choose 2$ pairs of
vertices, we have
$$ \sum_{i=1}^{n}({{d^-_i} \choose 2}+ d^-_i) \geq{ n \choose
2}.$$ Therefore, we have $\sum_{i=1}^{n} ((d^-_i)^2 + d^-_i) \geq
n(n-1)$. So, $n (d^2+d) \geq n(n-1)$ and $ d \geq
\frac{1}{2}(\sqrt{4n-3}-1)$ and the result follows by integrality
of $d$.\qed

\2

The digraph $H$ in Fig. \ref{f} shows that $\mu(6)=2.$ Indeed,
$f(6)=2=\Delta^-(H).$

\section{Families of Sets and
$\mu(n)$}\label{intersecting}

Since we will heavily use the terminology and  results of design
theory, in this section we characterize $\mu(n)$ in terms of
special families of sets. Symmetric families, 2-covering families
and families having a system of distinct representatives are of
significant interest in the theory and applications of
combinatorics, see, e.g., \cite{beth1986,cameron1994,jukna2001}.

We consider families of subsets of a finite set $X.$ Using
block-design terminology, we call the elements of $X$ {\em points}
and the subsets of $X$ {\em blocks}. Let ${\cal
F}=\{X_1,X_2,\ldots,X_m\}$ be a family. An $m$-tuple
$S=(x_1,x_2,\ldots,x_m)$ is a {\em system of distinct
representatives} ({\em SDR}) if all points of $S$ are distinct and
$x_i\in X_i$ for each $i=1,2,\ldots,m.$ A family $\cal F$ is {\em
symmetric} if $m=|X|.$ A family $\cal F$ is {\em 2-covering} if,
for each pair $j,k\in X,$ there exists $i\in \{1,2,\ldots,m\}$
such that $\{j,k\}\subseteq  X_i.$

Let $[n]=\{1,2,\ldots,n\}.$ Let ${\rm mcard}({\cal F})$ be the
maximum cardinality of a block in $\cal F.$ We call ${\cal F}$
{\em mediated} if it is symmetric, 2-covering and has an SDR. Let
$\mu^-(n)$ be the minimum ${\rm mcard}({\cal F})$ over all
mediated families on $[n].$

We have the following:

\begin{proposition}\label{muplus} For each $n\ge 1$, $\mu(n)=\mu^-(n)-1.$
\end{proposition}
\pf Let $D$ be a mediated digraph on vertices $[n]$ with
$\Delta^-(D)=\mu(n)$. By the definition of a mediated digraph, the
family ${\cal N}=\{N^-[i]:\ i\in [n]\}$ is 2-covering. Clearly,
$(1,2,\ldots,n)$ is an SDR of ${\cal N}.$ Thus, $\cal N$ is
mediated and $\mu^-(n)\le \mu(n)+1.$

Let ${\cal F}=\{X_1,X_2,\ldots,X_n\}$ be a mediated family on
$[n]$ with ${\rm mcard}({\cal F})=\mu^-(n).$ Since $\cal F$ has an
SDR (since it is mediated), without loss of generality, we may
assume that $i\in X_i.$ Construct a digraph $D$ with $V(D)=[n]$
and $N^-[i]=X_i.$ Since $\cal F$ is 2-covering, $D$ is mediated
and $\mu(n)\le \mu^-(n)-1.$ This inequality and $\mu^-(n)\le
\mu(n)+1$ imply that $\mu(n)=\mu^-(n)-1.$ \qed

\2

Let $n>k\ge 2$ and $\lambda \ge 1$ be integers. A family ${\cal
F}=\{X_1,X_2,\ldots,X_b\}$ of blocks on $X$ is called an
$(n,k,\lambda)$-{\em design} if $|X|=n$, each block has $k$ points
and every pair of distinct points is contained in exactly
$\lambda$ blocks. An $(n,k,\lambda)$-design is {\em symmetric} if
it has $n$ blocks, i.e., $b=n$. A {\em projective plane of order}
$q$ is a symmetric $(q^2+q+1,q+1,1)$-design for some integer
$q>1.$ For a family $\cal F$ of blocks and a point $i$, let $d(i)$
denote the number of blocks containing $i.$

The following two theorems are well-known, see, e.g.,
\cite{beth1986,cameron1994,jukna2001}.

\begin{theorem}\label{existpp}
For each prime power $q$, there exists a projective plane of order
$q$.
\end{theorem}

\begin{theorem}\label{regular}
Let ${\cal S}=\{X_1,X_2,\ldots,X_n\}$ be a family of subsets of
$\{1,2,\ldots,n\}$ and let $r$ be a natural number such that
$|X_i|=d(i)=r$ for each $i=1,2,\ldots,n$. Then $\cal S$ has an
SDR.
\end{theorem}

The last theorem can be used to prove the following:

\begin{proposition}\label{designcom}
Every symmetric $(n,k,\lambda)$-design is a mediated family of
blocks.
\end{proposition}
\pf Let ${\cal F}=\{X_1,X_2,\ldots,X_b\}$ be an
$(n,k,\lambda)$-design on $X$, $|X|=n$. It is well-known (see,
e.g., \cite{beth1986,cameron1994}) that, for all such designs,
there is a constant $r$ such that $r=d(i)$ for each point $i$. The
parameters $n,k,\lambda,b$ and $r$ also satisfy the following two
equalities: $bk(k-1)=\lambda n(n-1)$ and $r(k-1)=\lambda (b-1).$
Assume that $\cal F$ is symmetric. Using $b=n$ and the two
equalities, we easily conclude that $r=k.$ It now follows from
Theorem \ref{regular} that $\cal F$ has a SDR. Since $\cal F$ is
symmetric and 2-covering ($\lambda \ge 1$), $\cal F$ is
mediated.\qed

\2

Now we are ready to compute an infinite number of values of
$\mu(n).$

\begin{corollary}\label{muproj}
For each prime power $q$, $\mu(q^2+q+1)=f(q^2+q+1)=q.$
\end{corollary}
\pf Let $n=q^2+q+1.$ By Theorem \ref{existpp} and Propositions
\ref{lbprop} and \ref{designcom}, we have $f(n)\le
\mu(n)=\mu^-(n)-1\le q.$ However, one can trivially verify that
$f(n)=q.$\qed

\section{Upper Bounds for $\mu(n)$}\label{upper}

\begin{theorem}\label{ub}
Let $n=q^2+q+1+m(q+1)-t$, where $q$ is a prime power, $1 \leq m
\leq q+1$ and $0 \leq t \leq q$. Then $\mu(n) \leq q+m$.
\end{theorem}
\pf By Theorem \ref{existpp}, there exists a projective plane,
$\Pi$, of order $q$. Since $\Pi$ is a symmetric
$(q^2+q+1,q+1,1)$-design, $\Pi$ has $q^2+q+1$ blocks and $q^2+q+1$
points, each block has $q+1$ points and every point is contained
in $q+1$ blocks.

Let $P$ be the set of points in $\Pi$, let $x$ be a point in $\Pi$
and let $B_1,B_2, \ldots,B_{q+1}$ be the blocks of $\Pi$ which
contain $x$.

Let $W=\{ w_1,w_2, \ldots ,w_m\}$ be a set of extra points outside
the plane $\Pi.$ Let $Z=\{z_1,z_2,\ldots,z_{mq-t}\}$ be a subset
of $B_1 \cup B_2 \cup \ldots \cup B_m - \{ x \}$. Let
$Z'=\{z'_1,z'_2,\ldots,z'_{mq-t}\}$ such that $(P\cup W)\cap
Z'=\emptyset$. Since $\Pi$ is a design with $\lambda=1$, a point
$z$ in $Z$ is contained in exactly one of the blocks
$B_1,B_2,\ldots,B_{q+1}$, which we denote $B_{\tau(z)}$. We will
now define $n$ new blocks on the set $S=P \cup W \cup Z'$, such
that every pair of points in $S$ belong to a new block and no
block contains more than $q+1+m$ points. We will also see that the
new family of blocks has an SDR.

First add $W$ to all blocks $B_1,B_2, \ldots,B_{q+1}$. Then add
$z'$ to all blocks of $\Pi$ which contain $z$ but not $x$, for all
$z \in Z$. We include all extended and unextended blocks from
$\Pi$ to our new family of blocks. We thus have a set of $q^2+q+1$
blocks. We now add the following $m+|Z|$ blocks.

\2

 $Q_i = W \cup \{z'\in Z':\ z\in Z \cap B_i\}$, for $i=1,2,\ldots , m$.

 For every $z \in Z$ let $R_z = (B_{\tau(z)} - \{z\}) \cup \{z'\}.$

\2

By Proposition \ref{designcom}, $\Pi$ has an SDR. Thus, the new
blocks apart from the last $m+|Z|$ blocks have an SDR consisting
of points from P. This SDR can be extended to an SDR for all new
blocks by adding points from $W$ for blocks $Q_i$ and $Z'$ for
blocks $R_z.$

We will consider all possible pairs $\al,\be$ in $S$ and show that
for each there is a block containing $\al$ and $\be.$ This will
prove that the new family of blocks that we have constructed is
2-covering. We consider all possible cases for $\al,\be$ as follows.

{\bf Case 1:} $\al=a \in P$.

(i) If $\be=b \in P$, then some block contains both $a$ and $b$, as
$a,b\in P$, all pairs in $P$ are in a block of $\Pi$, and we have
either kept untouched or extended the blocks of $\Pi.$

(ii) If $\be=b' \in Z'-a'$ (if $a\not\in Z$, $a'=\emptyset$), then
$a$ and $b'$ both lie in some block, because of the following
argument. Some block must contain both $a$ and $b$, and if this
block does not contain $x$, then we have added $b'$ to this block,
and if it does contain $x$, then $a$ and $b'$ both lie in the blocks
$R_b$.

(iii) If $\be=b'=a'$, then $a$ and $b'$ both lie in all blocks
containing $a$ except the $B_i\cup W$'s.

(iv) If $\be=w \in W$, then $a$ and $w$ both lie in some block $B_i
\cup W$ since the sets of points $B_1,B_2,\ldots,B_{q+1}$ include
all points in $\Pi$ (each of these sets has $q+1$ points, and the
unique common point $x$).

{\bf Case 2:} $\al=w \in W$.

(i) If $\be=b \in W$, then $w$ and $b$ both lie in all blocks $B_i
\cup W$.

(ii) If $\be=b' \in Z'$, then $w$ and $b'$ both lie in some block
$Q_i$.

{\bf Case 3:} $\al=a' \in Z'$ and $\be=b' \in Z'$. Then $a'$ and
$b'$ both lie in some block, because of the following argument. Some
block must contain both $a$ and $b$, and if this block does not
contain $x$, then we have added both $a'$ and $b'$ to this block,
and if it does contain $x$, then $a'$ and $b'$ both lie in one of
the blocks $Q_i$.

To complete the proof it suffices to show that no new block has
size greater than $q+1+m$. (This puts an upper bound on the
maximum cardinality of the blocks and so an upper bound on
$\mu^-(n)$.) This is clearly true for all $Q_i$, $R_z$ and all
$B_i \cup W$. Now the proof follows because no block of $\Pi$ not
in the set $\{B_1,B_2,\ldots,B_{q+1}\}$ contains more than $m$
points from $B_1 \cup B_2 \cup \ldots \cup B_{m}$.\qed

\begin{corollary}\label{projth}
Let $q$ be a prime power. If $s$ is an integer such that
$q^2+q+2\le s\le q^2+2q+2,$ then $\mu(s)=f(s)=q+1.$
\end{corollary}
\pf Let $s$ be an integer such that $q^2+q+2\le s\le q^2+2q+2.$ By
Theorem \ref{ub} for $m=1$, $\mu(s)\le q+1.$ By Proposition
\ref{lbprop}, $q+1\ge \mu(s)\ge f(s)\ge f(q^2+q+2).$ Thus, it
suffices to show that $f(q^2+q+2)=q+1$, which is easily
verifiable. \qed

The following number-theoretical result was proved in
\cite{bakerPLMS83}.

\begin{theorem}\label{gapl}
There is a real $x_0$ such that for all $x>x_0$ the interval
$[x,x+x^{\alpha}]$, where $\alpha = 0.525$, contains prime
numbers.\end{theorem}

The last two assertions imply the following:

\begin{theorem}\label{o(1)th}
We have $\mu(n)=f(n)(1+o(1)).$
\end{theorem}
\pf Let $n$ be sufficiently large. Let $p$ and $q$ be a pair of
consecutive primes such that $p^2+p+1\le n<q^2+q+1$, and let
$d=q^2+q-p^2-p.$ By Theorem \ref{ub}, $\mu(n)\le p+\lceil d/(p+1)
\rceil.$ By Theorem \ref{gapl}, $q-p\le p^{\alpha}$. Thus,
$d=(q+p+1)(q-p)\le 3p\times p^{\alpha}=3p^{1+\alpha}.$ So,
$\mu(n)\le p + 3p^{\al}+1=p(1+o(1))=f(p^2+p+1)(1+o(1))\le
f(n)(1+o(1)).$\qed

\2

We believe that the following holds for a small constant $c$:

\begin{conjecture}\label{cconj}
There is a constant $c$ such that $\mu(n)\le f(n)+c$ for each $n.$
\end{conjecture}

If this conjecture holds, we would like to know the smallest value
of $c.$

To obtain another upper bound for $\mu(n)$ we will use the notion
of a cyclic $n$-difference cover that extends that of a cyclic
$(n,k,\lambda)$-difference set (see \cite{beth1986,jukna2001}). A
subset $D=\{d_1,d_2,\ldots,d_k\}$ of $\mathbb{Z}_n$ is called a
{\em cyclic} $n$-{\em  difference cover} if the collection of
values $d_i-d_j$ (mod $n$) contains every element of $Z_n$ at
least once. In the rest of this section, all operations with
elements of $\mathbb{Z}_n$ are taken modulo $n.$ For $c\in
\mathbb{Z}_n$, let $c+D=\{c+d:\ d\in D\}.$ The family ${\rm dev}
D=\{c+D:\ c\in \mathbb{Z}_n\}$ of $n$ blocks is called the {\em
development} of $D$.

\begin{proposition}\label{dcover}
If there exists a cyclic $n$-difference cover
$D=\{d_1,d_2,\ldots,d_k\}$, then $\mu(n)\le k-1.$
\end{proposition}
\pf Let $D=\{d_1,d_2,\ldots,d_k\}$ be a cyclic $n$-difference
cover. Consider ${\rm dev} D.$ Clearly, ${\rm dev} D$ is symmetric
and has an SDR $(d_1,d_1+1,\ldots,n-1+d_1)$.

For an arbitrary pair $a,b$ of distinct elements in
$\mathbb{Z}_n$, $a-b\in \mathbb{Z}_n$. Thus, there are $d_i,d_j\in
D$ such that $d_i-d_j=a-b$. Let $a=c+d_i,b=c'+d_j,$ where $c,c'\in
\mathbb{Z}_n$. The last three equalities imply $c=c'.$ Therefore,
$a$ and $b$ are both in $c+D.$ Hence, ${\rm dev} D$ is 2-covering
and, thus, mediated. So, by Proposition \ref{muplus}, $\mu(n)\le
k-1.$\qed

\2

Using a computer search the authors of
\cite{haanpaaJIS7,wiedemannCN90} determined the least $k=k(n)$
such that there is a cyclic $n$-difference cover
$\{d_1,d_2,\ldots,d_k\}$ for each $n\in \{3,4,5,\ldots,133\}$.
These results show that $\mu(n)-f(n)\le 1$ for $n\le 133,$ which
provides some support to Conjecture \ref{cconj}. However, for
large values of $n$, the bound in Proposition \ref{dcover} may not
be of much value since no upper bound on $k(n)$ of the form
$\sqrt{n}(1+o(1))$ seems to be known (see \cite{colbournIPL75}
where the bound $k(n)\le \sqrt{1.5n}+6$ was proved) and perhaps
the bound $k(n)\le \sqrt{n}(1+o(1))$ is simply not true.

\section{When $\mu(n)>f(n)$}\label{>sec}

Corollaries \ref{muproj} and  \ref{projth} may prompt some to
suspect that $\mu(n)=f(n)$ holds for each $n\ge 1.$ However, this
is not the case.

One of the best known conjectures in combinatorics is that a
projective plane does not exist if $q$ is not a prime power. The
celebrated Bruck-Ryser theorem \cite{bruckCJM1} (see also, e.g.,
\cite{cameron1994}) proves that if a projective plane of order $q$
exists, where $q\equiv 1$ or 2 (mod 4), then $q$ is the sum of two
squares of integers. This gives infinitely many values of $q$ for
which there is no projective plane of order $q$ (for example,
every number $q=2p,$ where $p$ is a prime congruent to 3 mod 4).
The fact that there are infinitely many primes congruent to 3 mod
4 follows from the famous Dirichlet's theorem: every arithmetic
progression with common difference relatively prime to the initial
term contains infinitely many prime numbers (see, e.g.,
\cite{narkiewicz1983}). The above implies the following:

\begin{theorem}\label{nopp}
There are infinitely  many positive integers $q$ for which there
is no projective plane of order $q.$
\end{theorem}

The non-existence of a projective plane of order 10, which does
not follow from the Bruck-Ryser theorem, was proved in
\cite{lamCJM41}.

\begin{theorem} \label{mu>f}
If there is no projective plane of order $q$, then $\mu(q^2+q+1) >
f(q^2+q+1)$.
\end{theorem}
\pf Let $q$ be an integer such that there is no projective plane
of order $q$, and let $n=q^2+q+1$. Suppose that $\mu(n) = f(n)$.
Observe that $f(n)=q.$ Thus, by Proposition \ref{muplus},
$\mu^-(n)=f(n)+1=q+1.$ Let ${\cal L}=\{L_1,L_2,\ldots ,L_n\}$ be a
mediated family of subsets of $[n]=\{1,2,\ldots,n\}$ with ${\rm
mcard}({\cal L})=q+1.$

We will obtain a contradiction by showing that ${\cal L}$ must be
a projective plane. By the choice of $n$, it suffices to prove
that $|L_i \cap L_j|=1$ for all $1 \leq i < j \leq n$ and
$|L_i|=q+1$ for each $i\in [n].$

Define $Q$ as follows, $Q=\{ \{i,j,L_k\}:\ \{i,j\} \subseteq L_k,
k\in [n]\}$. Observe that $L_k$ contains $|L_k|(|L_k|-1)/2$ pairs
of distinct points, so $|Q|= \sum^n_{k = 1} |L_k|(|L_k|-1)/2$.
Since $\cal L$ is mediated, every pair of points $i,j$ will appear
at least once in $Q$, so $|Q| \geq n(n-1)/2$. As $|L_k| \leq q+1$
for every $k \in [n]$, we have the following:

$$n\frac{(q+1)q}{2} \geq \sum_{k=1}^n \frac{|L_k|(|L_k|-1)}{2}= |Q| \geq \frac{n(n-1)}{2} \\
= n \frac{(q+1)q}{2}.$$

This implies that we must have equality everywhere, and thus
$|L_k|=q+1$ for each $k\in [n]$ and $|L_i \cap L_j|=1$ for all $1
\leq i < j \leq n$. \qed

This theorem and Theorem \ref{nopp} imply the following:

\begin{corollary}
For an infinite number of values of $n$, $\mu(n) > f(n).$
\end{corollary}

The following  problem is of certain interest.
\begin{problem}
Is $\mu(n)\le \mu(n+1)$ for each $n$?
\end{problem}

\2

{\bf Acknowledgements.} We'd like to thank Peter Cameron,
Christian Elsholtz, Adrian Sanders and Vsevolod Lev for
discussions on the topic of the paper.

\end{document}